\documentclass[11pt,a4paper,leqno]{amsart}
\usepackage{hyperref,cite}
\usepackage{graphics}
\usepackage{amsmath}
\newtheorem{thmx}{Theorem}[section]
\newtheorem{theorem}{Theorem}[section]
\newtheorem{remark}[theorem]{Remark}
\newtheorem{corollary}[theorem]{Corollary}
\theoremstyle{theorem}

\newtheorem{lemma}[theorem]{Lemma}
\theoremstyle{definition}

\numberwithin{equation}{section} \makeatletter
\begin{document}

\title{Some extensions of Enestr\"{o}m-Kakeya Theorem}

\maketitle

\begin{center}
N.A. Rather$^1$, Suhail Gulzar$^{2,*}$ and S.H. Ahangar$^2$
\end{center}
\begin{center}
$^1$Department of Mathematics,\\
University of Kashmir,\\
Hazratbal Srinagar 190006, India
\end{center}
\begin{center}
$^2$Department of Computer Science \& Engineering,\\
Islamic University of Science \& Technology,\\
Awantipora, Kashmir 192122, India\\
$^*$email:sgmattoo@gmail.com 
\end{center}
\begin{abstract}
 In this paper we obtain some refinements of a well-known result of Enestr\"{o}-Kakeya concerning the bounds for the moduli of the zeros of polynomials with complex
coefficients which improve upon some results due to Aziz and Mohammad, Govil and Rahman and others.
\end{abstract}
\section{Introduction}
A classical result due to Enestr\"{o}m (see\cite{e}) and Kakeya \cite{sk} concerning the bounds for the moduli of the zeros of polynomials having postive and real coefficients is often stated as in the following theorem (see \cite[p.136]{marden} or \cite[p.272]{mmr}).
\begin{thmx}
If $P(z)=\sum_{j=0}^{n}a_{j}z^{j}$ is a polynomial of degree $ n $ with real coefficients satisfying
$$ a_{n}\geq a_{n-1}\geq \cdots \geq a_{1}\geq a_{0}>0,$$
then all its zeros lie in $ |z|\leq 1.$
\end{thmx}
In literature \cite{asv79,asv81,cs,gr,km,k,marden,mmr,rsm,rsuh} there exits several extensions of Enestr\"{o}m-Kakeya theorem. 
A. Aziz and Q. G. Muhammad \cite{am} used Schwarz's Lemma and proved among other things the following generalization of Enestr\"{o}m-Kakeya Theorem. 
\begin{thmx}\label{ta}
Let $P(z)=\sum_{j=0}^{n}a_jz^n$ be a polynomial of degree $n$ with real coefficients. If $t_1>t_2\geq 0$ can be found such that
$$ t_1t_2a_r+(t_1-t_2)a_{r-1}-a_{r-2}\geq 0\quad\textnormal{for}\quad r=1,2,\cdots,n+1\quad(a_{-1}=a_{n+1}=0),  $$
then all the zeros of $P(z)$ lie in $|z|\leq t_1.$
\end{thmx}
N.K. Govil and Q. I. Rahman \cite{gr} extended the Enestr\"{o}m-Kakeya theorem to the polynomials with complex coefficients by proving:
\begin{thmx}\label{tgr}
Let $P(z)=\sum_{j=0}^{n}a_jz^n$ be a polynomial of degree $n$ with complex coefficients such that
$$|\arg a_\nu-\beta|\leq \alpha\leq \pi/2,\quad \nu=0,1,2,\cdots,n,   $$
for some real $\beta$ and
$$|a_n|\geq |a_{n-1}|\geq \cdots|a_1|\geq |a_0|,  $$
then all the zeros of $P(z)$ lie in 
$$|z|\leq (\cos\alpha+\sin\alpha)+\dfrac{2\sin\alpha}{|a_n|}\sum\limits_{\nu=0}^{n}|a_\nu|.  $$
\end{thmx}
The following generalization of Theorem 1.3 is due to A. Aziz and Q. G. Mohammad \cite{am},
\begin{thmx}\label{tc'}
Let $P(z)=\sum_{j=0}^{n}a_jz^n$ be a polynomial of degree $n$ with real coefficients. If $t_1>t_2\geq 0$ can be found such that
$$|\arg a_\nu-\beta|\leq \alpha\leq \pi/2,\qquad\qquad\nu=0,1,2,\cdots,n  $$ 
for some real $\beta$ and for some $t>0,$
\begin{align*}
 t^n|a_n|\leq &t^{n-1}|a_{n-1}|\leq\cdots\leq t^{k+1}|a_{k+1}|\leq t^{k}|a_k|,\\&t^{k}|a_k|\geq t^{k-1}|a_{k-1}|\geq\cdots\geq t|a_1|\geq |a_0| >0, 
 \end{align*}
 where $0\leq k\leq n,$ then $P(z)$ has all its zeros in the circle
 $$|z|\leq t\left\{ \left(\dfrac{2t^k|a_k|}{t^n|a_n|}-1\right)\cos\alpha+\sin\alpha\right\}+2\sin\alpha\sum\limits_{j=0}^{n}\dfrac{|a_j|}{|a_n|t^{n-j-1}}.  $$
\end{thmx}
Rather et al. \cite{rsm} extended Theorem 1.2 to the polynomials with complex coefficients and proved the following result.
\begin{thmx}\label{tb}
Let $P(z)=\sum_{j=0}^{n}a_jz^n$ be a polynomial of degree $n$ with complex coefficients such that
$$|\arg a_\nu-\beta|\leq \alpha\leq \pi/2,\quad \nu=0,1,2,\cdots,n,   $$
for some real $\beta$. If $t_1>t_2\geq 0$ can be found such that
$$ t_1t_2|a_r|+(t_1-t_2)|a_{r-1}|-|a_{r-2}|\geq 0\quad\textnormal{for}\quad r=1,2,\cdots,k+1,  $$
and
$$ t_1t_2|a_r|+(t_1-t_2)|a_{r-1}|-|a_{r-2}|\leq 0\quad\textnormal{for}\quad r=k+2,\cdots,n+1,  $$
$0\leq k\leq n,$ $a_{-1}=a_{n+1}=0,$ then all the zeros of $P(z)$ lie in 
\begin{align}\label{tbe}
|z|\leq t_1\left\{\left(\dfrac{2|a_k|+2t_2|a_{k+1}|}{t_1^{n-k}|a_n|}-1\right)\cos\alpha+\sin\alpha\right\}+2\dfrac{\sin\alpha}{|a_n|}\sum\limits_{j=0}^{n-1}\dfrac{|a_j|}{t_1^{n-j-1}}.
\end{align}
\end{thmx}
They \cite{rsm} also obtained the following generalization of Theorem 1.2.
\begin{thmx}\label{tc}
Let $P(z)=\sum_{j=0}^{n}a_jz^n$ be a polynomial of degree $n$ with $\textnormal{Re}\, a_j=\alpha_j$ and $\textnormal{Im}\, a_j=\beta_j,$ $j=0,1,\cdots,n.$ If $t_1>t_2\geq 0$ can be found such that
$$ t_1t_2\alpha_r+(t_1-t_2)\alpha_{r-1}-\alpha_{r-2}\geq 0\quad\textnormal{for}\quad r=1,2,\cdots,k+1,  $$
$$ t_1t_2\alpha_r+(t_1-t_2)\alpha_{r-1}-\alpha_{r-2}\leq 0\quad\textnormal{for}\quad r=k+2,\cdots,n+1,  $$
and
$$ t_1t_2\beta_r+(t_1-t_2)\beta_{r-1}-\beta_{r-2}\geq 0\quad\textnormal{for}\quad r=1,2,\cdots,m+1,  $$
$$ t_1t_2\beta_r+(t_1-t_2)\beta_{r-1}-\beta_{r-2}\leq 0\quad\textnormal{for}\quad r=m+2,\cdots,n+1,  $$
$0\leq k\leq n,$ $0\leq m\leq n,$ $\alpha_{-1}=\beta_{-1}=\alpha_{n+1}=\beta_{n+1}=0,$ $\alpha_n>0,$ then all the zero of $P(z)$ lie in 
\begin{align}\label{tec}
|z|\leq \dfrac{t_1}{|a_n|}\left\{2t_1^{k-n}\big(\alpha_{k}+t_2\alpha_{k+1}\big)+2t_1^{m-n}\big(\beta_{m}+t_2\beta_{m+1}\big)-\big(\alpha_n+\beta_n\big) \right\}. 
\end{align}
\end{thmx}
In this paper, we shall first present the following result which among other things considerably improves the bound of Theorem 1.5 for $0\leq k\leq n-3.$
\begin{theorem}\label{t1}
Let $P(z)=\sum_{j=0}^{n}a_jz^n$ be a polynomial of degree $n\geq 3$ with complex coefficients such that
$$|\arg a_j-\beta|\leq \alpha\leq \pi/2\quad\textnormal{for},\quad j=0,1,2,\cdots,n   $$
for some real $\beta.$ If $t_1(\neq 0)$ and $t_2$ with $t_1\geq t_2\geq 0$ can be found such that
$$ t_1t_2|a_r|+(t_1-t_2)|a_{r-1}|-|a_{r-2}|\geq 0\quad\textnormal{for}\quad r=1,2,\cdots,k+1,  $$
and
$$ t_1t_2|a_r|+(t_1-t_2)|a_{r-1}|-|a_{r-2}|\leq 0\quad\textnormal{for}\quad r=k+2,\cdots,n+1,  $$
$0\leq k\leq n-3,$ $a_{-1}=a_{n+1}=0,$ then all the zeros of $P(z)$ lie in 
\begin{align}\label{t1e}\nonumber
\left|z+\dfrac{a_{n-1}}{a_n}-(t_1-t_2)\right|\leq&\,\,\dfrac{2t_2|a_{k+1}|+2|a_k|}{|a_n|t_1^{n-k-1}}\cos\alpha+2\dfrac{\sin\alpha}{|a_n|}\sum\limits_{\nu=0}^{n-2}\dfrac{|a_{\nu}|}{t_1^{n-\nu-1}}\\&+\left(t_2+\left|\dfrac{a_{n-1}}{a_n}\right|\right)(\sin\alpha-\cos\alpha).
\end{align}
\end{theorem}

\begin{remark}
\textnormal{
In general Theorem 1.7 gives much better result than Theorem 1.5 for $0\leq k\leq n-3.$ To see this, we show that the circle defined by \eqref{t1e} is contained in the circle defined by \eqref{tbe}. Let $z=w$ be any point belonging to the circle defined by \eqref{t1e}, then
\begin{align*}\nonumber
\left|w+\dfrac{a_{n-1}}{a_n}-(t_1-t_2)\right|\leq&\left(t_2+\left|\dfrac{a_{n-1}}{a_n}\right|\right)(\sin\alpha-\cos\alpha)\\&+\dfrac{2t_2|a_{k+1}|+2|a_k|}{|a_n|t_1^{n-k-1}}\cos\alpha+2\dfrac{\sin\alpha}{|a_n|}\sum\limits_{\nu=0}^{n-2}\dfrac{|a_{\nu}|}{t_1^{n-\nu-1}}.
\end{align*}
This implies,
\begin{align*}\nonumber
|w|=&\bigg|w+\dfrac{a_{n-1}}{a_n}-(t_1-t_2)+(t_1-t_2)-\dfrac{a_{n-1}}{a_n}\bigg|\\\nonumber\leq& \bigg|w+\dfrac{a_{n-1}}{a_n}-(t_1-t_2)\bigg|+\bigg|(t_1-t_2)-\dfrac{a_{n-1}}{a_n}\bigg|\\\leq &\left(t_2+\left|\dfrac{a_{n-1}}{a_n}\right|\right)(\sin\alpha-\cos\alpha)+\dfrac{2t_2|a_{k+1}|+2|a_k|}{|a_n|t_1^{n-k-1}}\cos\alpha\\&+2\dfrac{\sin\alpha}{|a_n|}\sum\limits_{\nu=0}^{n-2}\dfrac{|a_{\nu}|}{t_1^{n-\nu-1}}+\bigg|(t_1-t_2)-\dfrac{a_{n-1}}{a_n}\bigg|.
\end{align*}
Using Lemma 2.1 with $j=n+1$ and noting that $a_{n+1}=0,$ we get
\begin{align*}
|w|\leq&\left(t_2+\left|\dfrac{a_{n-1}}{a_n}\right|\right)(\sin\alpha-\cos\alpha)+\dfrac{2t_2|a_{k+1}|+2|a_k|}{|a_n|t_1^{n-k-1}}\cos\alpha\\&\qquad+2\dfrac{\sin\alpha}{|a_n|}\sum\limits_{\nu=0}^{n-2}\dfrac{|a_{\nu}|}{t_1^{n-\nu-1}}\\\nonumber&+\left(\left|\dfrac{a_{n-1}}{a_n}\right|-(t_1-t_2)\right)\cos\alpha+\left(t_1-t_2+\left|\dfrac{a_{n-1}}{a_n}\right|\right)\sin\alpha\\\nonumber=&2\left(\dfrac{t_2|a_{k+1}|+|a_k|}{t_1^{n-k-1}|a_n|}\right)\cos\alpha+2\dfrac{\sin\alpha}{|a_n|}\sum\limits_{\nu=0}^{n-2}\dfrac{|a_{\nu}|}{t_1^{n-\nu-1}}\\&\qquad+2\sin\alpha\dfrac{|a_{n-1}|}{|a_n|}+t_1(\sin\alpha-\cos\alpha)\\& =t_1\left\{\left(\dfrac{2|a_k|+2t_2|a_{k+1}|}{t_1^{n-k}|a_n|}-1\right)\cos\alpha+\sin\alpha\right\}+2\dfrac{\sin\alpha}{|a_n|}\sum\limits_{j=0}^{n-1}\dfrac{|a_j|}{t_1^{n-j-1}}.
\end{align*}
This shows that the point $z=w$ also belongs to the circle defined by \eqref{tbe}. Hence the circle defined by \eqref{t1e} is contained in the circle defined by \eqref{tbe}.
}
\end{remark}

For $t_2=0,$ in Theorem 1.7, we obtain the following result which considerably improves the bound of Theorem 1.4 for $0\leq k\leq n-3.$
\begin{corollary}
Let $P(z)=\sum_{j=0}^{n}a_jz^n$ be a polynomial of degree $n\geq 3$ with complex coefficients such that
$$|\arg a_j-\beta|\leq \alpha\leq \pi/2\quad\textnormal{for},\quad j=0,1,2,\cdots,n   $$
for some real $\beta.$ If $t> 0$ can be found such that
\begin{align*}
 t^n|a_n|\leq &t^{n-1}|a_{n-1}|\leq\cdots\leq t^{k+1}|a_{k+1}|\leq t^{k}|a_k|,\\&t^{k}|a_k|\geq t^{k-1}|a_{k-1}|\geq\cdots\geq t|a_1|\geq |a_0| >0 
 \end{align*}
$0\leq k\leq n-3,$ $a_{-1}=a_{n+1}=0,$ then all the zeros of $P(z)$ lie in 
\begin{align}\nonumber
\left|z+\dfrac{a_{n-1}}{a_n}-t\right|\leq&\left|\dfrac{a_{n-1}}{a_n}\right|(\sin\alpha-\cos\alpha)\\&+\dfrac{2|a_k|}{|a_n|t_1^{n-k-1}}\cos\alpha+2\dfrac{\sin\alpha}{|a_n|}\sum\limits_{\nu=0}^{n-2}\dfrac{|a_{\nu}|}{t_1^{n-\nu-1}}.
\end{align}
\end{corollary}
Next, as a generalization of Theorem 1.6, we prove the following result.
\begin{theorem}\label{t2}
Let $P(z)=\sum_{j=0}^{n}a_jz^n$ be a polynomial of degree $n\geq 1$ with $\textnormal{Re}\, a_j=\alpha_j$ and $\textnormal{Im}\, a_j=\beta_j,$ $j=0,1,\cdots,n.$ If $t_1\neq 0,$ $t_1\geq t_2\geq 0$ can be found such that
$$ t_1t_2\alpha_r+(t_1-t_2)\alpha_{r-1}-\alpha_{r-2}\geq 0\quad\textnormal{for}\quad r=1,2,\cdots,k+1,  $$
$$ t_1t_2\alpha_r+(t_1-t_2)\alpha_{r-1}-\alpha_{r-2}\leq 0\quad\textnormal{for}\quad r=k+2,\cdots,n+1,  $$
and
$$ t_1t_2\beta_r+(t_1-t_2)\beta_{r-1}-\beta_{r-2}\geq 0\quad\textnormal{for}\quad r=1,2,\cdots,m+1,  $$
$$ t_1t_2\beta_r+(t_1-t_2)\beta_{r-1}-\beta_{r-2}\leq 0\quad\textnormal{for}\quad r=m+2,\cdots,n+1,  $$ 
$0\leq k\leq n-1,$ $0\leq m-1\leq n,$ $\alpha_{-1}=\beta_{-1}=\alpha_{n+1}=\beta_{n+1}=0,$ $\alpha_n>0,$ then all the zeros of $P(z)$ lie in 
\begin{align}\label{t2e}\nonumber
\bigg|z+&\dfrac{\alpha_{n-1}-(t_1-t_2)\alpha_n}{a_n}\bigg|\\\nonumber&\leq 2(\alpha_{k+1}t_2+\alpha_k)\dfrac{t_1^{k+1}}{|a_n|t_1^{n}}+2(\beta_{m+1}t_2+\beta_m)\dfrac{t_1^{m+1}}{|a_n|t_1^{n}}\\&\qquad\qquad\qquad-\dfrac{t_2\alpha_n+t_1\beta_n+\alpha_{n-1}}{|a_n|}.
\end{align}
\end{theorem}
\begin{remark}
\textnormal{
In general Theorem 1.10 also gives much better result than Theorem 1.6 for $0\leq k\leq n-1.$ For this, we show that the circle defined by \eqref{t2e} is contained in the circle defined by \eqref{tec}. Let $z=w$ be any point belonging to the circle defined by \eqref{t2e} then
\begin{align*}\nonumber
\bigg|w+&\dfrac{\alpha_{n-1}-(t_1-t_2)\alpha_n}{a_n}\bigg|\\\nonumber&\leq 2(\alpha_{k+1}t_2+\alpha_k)\dfrac{t_1^{k+1}}{|a_n|t_1^{n}}+2(\beta_{m+1}t_2+\beta_m)\dfrac{t_1^{m+1}}{|a_n|t_1^{n}}\\&\qquad\qquad\qquad-\dfrac{t_2\alpha_n+t_1\beta_n+\alpha_{n-1}}{|a_n|}.
\end{align*}
This implies 
\begin{align*}
|w|=&\bigg|w+\dfrac{\alpha_{n-1}-(t_1-t_2)\alpha_n}{a_n}-\dfrac{\alpha_{n-1}-(t_1-t_2)\alpha_n}{a_n}\bigg|\\\leq&\bigg|w+\dfrac{\alpha_{n-1}-(t_1-t_2)\alpha_n}{a_n}\bigg|+\bigg|\dfrac{\alpha_{n-1}-(t_1-t_2)\alpha_n}{a_n}\bigg|\\\leq&2(\alpha_{k+1}t_2+\alpha_k)\dfrac{t_1^{k+1}}{|a_n|t_1^{n}}+2(\beta_{m+1}t_2+\beta_m)\dfrac{t_1^{m+1}}{|a_n|t_1^{n}}\\&\qquad\quad-\dfrac{t_2\alpha_n+t_1\beta_n+\alpha_{n-1}}{|a_n|}+\dfrac{\alpha_{n-1}-(t_1-t_2)\alpha_n}{|a_n|}\\=&\dfrac{t_1}{|a_n|}\left\{2t_1^{k-n}\big(\alpha_{k}+t_2\alpha_{k+1}\big)+2t_1^{m-n}\big(\beta_{m}+t_2\beta_{m+1}\big)-\big(\alpha_n+\beta_n\big) \right\}.
\end{align*}
Hence the point $z=w$ belongs to the circle defined by \eqref{tec} and therefore, the circle defined by \eqref{t2e} is contained in the circle defined \eqref{tec}.}
\end{remark}

For $\beta_j=0,j=0,1,\cdots,n$ in Theorem 1.10, we obtain the following result.
\begin{corollary}
Let $P(z)=\sum_{j=0}^{n}a_jz^n$ be a polynomial of degree $n\geq 1$ with real and positive coefficients.
 If $t_1(\neq 0)\geq t_2\geq 0$ can be found such that
$$ t_1t_2a_r+(t_1-t_2)a_{r-1}-a_{r-2}\geq 0\quad\textnormal{for}\quad r=1,2,\cdots,k+1,  $$
and
$$ t_1t_2a_r+(t_1-t_2)a_{r-1}-a_{r-2}\leq 0\quad\textnormal{for}\quad r=k+2,\cdots,n+1,  $$
$0\leq k\leq n-1,$ $a_{-1}=a_{n+1}=0,$ then all the zeros of $P(z)$ lie in 
\begin{align}\nonumber
\left|z+\dfrac{a_{n-1}}{a_n}-(t_1-t_2)\right|\leq&t_2+\dfrac{a_{n-1}}{a_n}+\dfrac{2t_2a_{k+1}+2a_k}{a_nt_1^{n-k-1}}
\end{align}
\end{corollary}
\vskip 6mm

\section{Preliminaries}
\begin{lemma}~\textnormal{(\cite{rsm})}\label{l2}
If $|\arg a_j-\beta|\leq \alpha\leq \pi/2,$ $j=0,1,2,\cdots,n$ and $\beta $ real, then for $t_1>t_2\geq 0,$
\begin{align*}
|t_1t_2a_j+(t_1-t_2)a_{j-1}&-a_{j-2}|\\\leq |t_1t_2&|a_j|+(t_1-t_2)|a_{j-1}|-|a_{j-2}||\cos\alpha\\\quad&+(t_1t_2|a_j|+(t_1-t_2)|a_{j-1}|+|a_{j-2}|)\sin\alpha.
\end{align*}
\end{lemma}
\section{Proof of theorems}
\begin{proof}[Proof of Theorem 1.7.]
Consider the polynomial
\begin{align}\nonumber
F(z)=&(t_1-z)(t_2+z)P(z)\\\nonumber=&-a_nz^{n+2}+(a_n(t_1-t_2)-a_{n-1})z^{n+1}\\&\nonumber\qquad+\sum_{\nu=2}^{n}(a_\nu t_1t_2+a_{\nu-1}(t_1-t_2)-a_{\nu-2})z^{\nu}\\\nonumber&\qquad\quad+(a_1t_1t_2+a_0(t_1-t_2))z+a_0t_1t_2\\\nonumber=&-a_nz^{n+2}+(a_n(t_1-t_2)-a_{n-1})z^{n+1}\\\nonumber\qquad\quad&+\sum_{\nu=0}^{n}(a_\nu t_1t_2+a_{\nu-1}(t_1-t_2)-a_{\nu-2})z^{\nu}\qquad (a_{-2}=a_{-1}=0)
\end{align}
Let $|z|>t_1,$ then
\begin{align}\nonumber\label{p1}
|F(z)|\geq&|a_n||z|^{n+1}\Bigg[\left|z+\dfrac{a_{n-1}}{a_n}-(t_1-t_2)\right|\\&\nonumber\quad-\sum_{\nu=0}^{n}\left|\dfrac{a_\nu t_1t_2+a_{\nu-1}(t_1-t_2)-a_{\nu-2}}{a_n}\right|\dfrac{1}{|z|^{n-\nu+1}}\Bigg]\\\nonumber>&|a_n||z|^{n+1}\Bigg[\left|z+\dfrac{a_{n-1}}{a_n}-(t_1-t_2)\right|\\&\quad-\dfrac{1}{|a_n|}\sum_{\nu=0}^{n}\dfrac{|a_\nu t_1t_2+a_{\nu-1}(t_1-t_2)-a_{\nu-2}|}{t_1^{n-\nu+1}}\Bigg]
\end{align}
By Lemma 2.1, we have
\begin{align*}
\sum_{\nu=0}^{n}&\dfrac{|a_\nu t_1t_2+a_{\nu-1}(t_1-t_2)-a_{\nu-2}|}{t_1^{n-\nu+1}}\\&\leq \sum_{\nu=0}^{n}\dfrac{||a_\nu| t_1t_2+|a_{\nu-1}|(t_1-t_2)-|a_{\nu-2}||}{t_1^{n-\nu+1}}\cos\alpha\\&\qquad+\sum_{\nu=0}^{n}\dfrac{||a_\nu| t_1t_2+|a_{\nu-1}|(t_1-t_2)+|a_{\nu-2}||}{t_1^{n-\nu+1}}\sin\alpha\\&=\sum_{\nu=0}^{k+1}\dfrac{|a_\nu| t_1t_2+|a_{\nu-1}|(t_1-t_2)-|a_{\nu-2}|}{t_1^{n-\nu+1}}\cos\alpha\\&\qquad-\sum_{\nu=k+2}^{n}\dfrac{|a_\nu| t_1t_2+|a_{\nu-1}|(t_1-t_2)-|a_{\nu-2}|}{t_1^{n-\nu+1}}\cos\alpha\\&\qquad+\sum_{\nu=0}^{n}\dfrac{|a_\nu| t_1t_2+|a_{\nu-1}|(t_1-t_2)+|a_{\nu-2}|}{t_1^{n-\nu+1}}\sin\alpha\\&=(t_2|a_n|+|a_{n-1}|)(\sin\alpha-\cos\alpha)\\&\qquad+\dfrac{2t_2|a_{k+1}|+2|a_k|}{t_1^{n-k-1}}\cos\alpha+2\sin\alpha\sum\limits_{\nu=0}^{n-2}\dfrac{|a_{\nu}|}{t_1^{n-\nu-1}}
\end{align*}
Using this in \eqref{p1}, we obtain
 \begin{align}\nonumber
 |F(z)|>&|a_n||z|^{n+1}\Bigg\{\left|z+\dfrac{a_{n-1}}{a_n}-(t_1-t_2)\right|\\\nonumber&-\left(t_2+\left|\dfrac{a_{n-1}}{a_n}\right|\right)(\sin\alpha-\cos\alpha)\\&\qquad-\dfrac{2t_2|a_{k+1}|+2|a_k|}{|a_n|t_1^{n-k-1}}\cos\alpha-2\sin\alpha\sum\limits_{\nu=0}^{n-2}\left|\dfrac{a_{\nu}}{a_n}\right|\dfrac{1}{t_1^{n-\nu-1}}\Bigg\}>0,
 \end{align}
whenever
\begin{align*}
\left|z+\dfrac{a_{n-1}}{a_n}-(t_1-t_2)\right|>&\left(t_2+\left|\dfrac{a_{n-1}}{a_n}\right|\right)(\sin\alpha-\cos\alpha)\\&+\dfrac{2t_2|a_{k+1}|+2|a_k|}{|a_n|t_1^{n-k-1}}\cos\alpha+2\sin\alpha\sum\limits_{\nu=0}^{n-2}\left|\dfrac{a_{\nu}}{a_n}\right|\dfrac{1}{t_1^{n-\nu-1}}.
\end{align*}
Hence, all the zeros of $F(z)$ whose modulus is greater than $t_1$ lie in the circle
\begin{align*}\nonumber
\left|z+\dfrac{a_{n-1}}{a_n}-(t_1-t_2)\right|\leq&\left(t_2+\left|\dfrac{a_{n-1}}{a_n}\right|\right)(\sin\alpha-\cos\alpha)\\&+\dfrac{2t_2|a_{k+1}|+2|a_k|}{|a_n|t_1^{n-k-1}}\cos\alpha+2\dfrac{\sin\alpha}{|a_n|}\sum\limits_{\nu=0}^{n-2}\dfrac{|a_{\nu}|}{t_1^{n-\nu-1}}.
\end{align*}
We now show that all those zeros of $F(z)$ whose modulus is less than or equal to $t_1$ also satisfy \eqref{t1e}. Let $|z|\leq t_1,$ then by using Lemma 2.1, we have
\begin{align*}
\left|z+\dfrac{a_{n-1}}{a_n}-(t_1-t_2)\right|&\leq t_1+\left|\dfrac{a_{n-1}}{a_n}-(t_1-t_2)\right|\\&\leq t_1+\left(\left|\dfrac{a_{n-1}}{a_n}\right|-(t_1-t_2)\right)\cos\alpha\\&\qquad\qquad\qquad+\left(t_1-t_2+\left|\dfrac{a_{n-1}}{a_n}\right|\right)\sin\alpha\\&\leq\left(t_2+\left|\dfrac{a_{n-1}}{a_n}\right|\right)(\sin\alpha-\cos\alpha)\\&\,\,+\dfrac{2t_2|a_{k+1}|+2|a_k|}{|a_n|t_1^{n-k-1}}\cos\alpha+2\sin\alpha\sum\limits_{\nu=0}^{n-2}\left|\dfrac{a_{\nu}}{a_n}\right|\dfrac{1}{t_1^{n-\nu-1}}
\end{align*}
if
\begin{align}\nonumber\label{p2}
t_1(1-\cos\alpha)&+2\left(t_2+\left|\dfrac{a_{n-1}}{a_n}\right|\right)\cos\alpha+t_1\sin\alpha\\&\leq \dfrac{2t_2|a_{k+1}|+2|a_k|}{|a_n|t_1^{n-k-1}}\cos\alpha+2\sin\alpha\sum\limits_{\nu=0}^{n-2}\left|\dfrac{a_\nu}{a_n}\right|\dfrac{1}{t_1^{n-\nu-1}}
\end{align}
Now, by given hypothesis,
\begin{align*}\nonumber
\sum_{\nu=k+2}^{n}\dfrac{|a_\nu| t_1t_2+|a_{\nu-1}|(t_1-t_2)-|a_{\nu-2}|}{t_1^{n-\nu+1}}\geq 0\qquad\qquad(0\leq k\leq n-2),
\end{align*}
that is,
\begin{align}\label{p3}
2\left(t_2+\left|\dfrac{a_{n-1}}{a_n}\right|\right)\leq \dfrac{2t_2|a_{k+1}|+2|a_k|}{|a_n|t_1^{n-k-1}}.
\end{align}
Therefore, \eqref{p2} holds, if
\begin{align}\label{p5}
t_1(1-\cos\alpha)+t_1\sin\alpha\leq 2\sin\alpha\sum\limits_{\nu=0}^{n-2}\left|\dfrac{a_\nu}{a_n}\right|\dfrac{1}{t_1^{n-\nu-1}}
\end{align}
Again, by hypothesis,
\begin{align}\label{p4}
 t_1\leq t_2+\frac{|a_{n-1}|}{|a_n|}
\end{align}
Using \eqref{p4} in \eqref{p3}, we obtain
\begin{align*}
t_1&\leq t_2+\frac{|a_{n-1}|}{|a_n|}\leq \dfrac{t_2|a_{k+1}|+|a_k|}{|a_n|t_1^{n-k-1}}\\& \leq \dfrac{t_1|a_{k+1}|+|a_k|}{|a_n|t_1^{n-k-1}}= \dfrac{|a_{k+1}|}{|a_n|t_1^{n-k-2}}+\dfrac{|a_{k}|}{|a_n|t_1^{n-k-1}}\\&\leq \sum\limits_{\nu=0}^{n-2}\left|\dfrac{a_\nu}{a_n}\right|\dfrac{1}{t_1^{n-\nu-1}}\qquad\qquad(0\leq k\leq n-3)
\end{align*}
Hence \eqref{p5} holds, if 
\begin{align*}\nonumber
t_1(1-\cos\alpha)+t_1\sin\alpha\leq 2t_1\sin\alpha,
\end{align*}
or if,
\begin{align}\label{p5'}
\cos\alpha+\sin\alpha\geq 1,\quad\textnormal{where}\quad 0\leq \alpha\leq \pi/2.
\end{align}
Note that when $0\leq \alpha\leq \pi/2,$
\begin{align*}
 \cos\alpha+\sin\alpha=\sqrt{2}\sin\left(\alpha+\pi/4\right)\geq \sqrt{2}.\dfrac{1}{\sqrt{2}}= 1
\end{align*}
This shows that \eqref{p5'} holds. Thus we have shown that if $|z|\leq t,$ then for $0\leq k\leq n-3$
\begin{align*}
\left|z+\dfrac{a_{n-1}}{a_n}-(t_1-t_2)\right|\leq&\left(t_2+\left|\dfrac{a_{n-1}}{a_n}\right|\right)(\sin\alpha-\cos\alpha)\\&+\dfrac{2t_2|a_{k+1}|+2|a_k|}{|a_n|t_1^{n-k-1}}\cos\alpha+2\dfrac{\sin\alpha}{|a_n|}\sum\limits_{\nu=0}^{n-2}\dfrac{|a_{\nu}|}{t_1^{n-\nu-1}}.
\end{align*}
Hence all the zeros of $F(z)$ lie in the circle defined by \eqref{t1e}. But all the zeros of $P(z)$ are also the zeros of $F(z),$ we conclude that all the zeros of $P(z)$ lie in the circle defined by \eqref{t1e}. This proves theorem 1.7. 
\end{proof}

\begin{proof}[Proof of Theorem 1.10.]
Consider the polynomial,
\begin{align}\nonumber
G(z)=&(t_1-z)(t_2+z)P(z)\\\nonumber=&-a_nz^{n+2}+(a_n(t_1-t_2)-a_{n-1})z^{n+1}\\&\nonumber\qquad+\sum_{\nu=2}^{n}(a_\nu t_1t_2+a_{\nu-1}(t_1-t_2)-a_{\nu-2})z^{\nu}\\\nonumber&\qquad\quad+(a_1t_1t_2+a_0(t_1-t_2))z+a_0t_1t_2\\\nonumber=&-a_nz^{n+2}+(a_n(t_1-t_2)-a_{n-1})z^{n+1}\\\nonumber\qquad\quad&+\sum_{\nu=0}^{n}(a_\nu t_1t_2+a_{\nu-1}(t_1-t_2)-a_{\nu-2})z^{\nu}\qquad (a_{-2}=a_{-1}=0)
\end{align}
Let $|z|>t_1,$ then
\begin{align}\nonumber\label{p31}
|G(z)|\geq&|z|^{n+1}\Bigg[\left|a_nz+a_{n-1}-(t_1-t_2)a_n\right|\\&\nonumber\quad-\sum_{\nu=0}^{n}\left|a_\nu t_1t_2+a_{\nu-1}(t_1-t_2)-a_{\nu-2}\right|\dfrac{1}{|z|^{n-\nu+1}}\Bigg]\\\nonumber>&|z|^{n+1}\Bigg[\left|a_nz+\alpha_{n-1}-(t_1-t_2)\alpha_n\right|-|\beta_{n-1}|-(t_1-t_2)|\beta_n|\\&\quad-\sum_{\nu=0}^{n}|a_\nu t_1t_2+a_{\nu-1}(t_1-t_2)-a_{\nu-2}|\dfrac{1}{t_1^{n-\nu+1}}\Bigg]
\end{align}
Now by hypothesis, 
\begin{align*}
\sum_{\nu=0}^{n}&|a_\nu t_1t_2+a_{\nu-1}(t_1-t_2)-a_{\nu-2}|t_{1}^{\nu}\\&\leq \sum_{\nu=0}^{n}|\alpha_\nu t_1t_2+\alpha_{\nu-1}(t_1-t_2)-\alpha_{\nu-2}|t_{1}^{\nu}\\&\qquad\qquad+\sum_{\nu=0}^{n}|\beta_\nu t_1t_2+\beta_{\nu-1}(t_1-t_2)-\beta_{\nu-2}|t_{1}^{\nu}\\&\leq\sum_{\nu=0}^{k+1}|\alpha_\nu t_1t_2+\alpha_{\nu-1}(t_1-t_2)-\alpha_{\nu-2}|t_{1}^{\nu}\\&+\sum_{\nu=k+2}^{n}|\alpha_\nu t_1t_2+\alpha_{\nu-1}(t_1-t_2)-\alpha_{\nu-2}|t_{1}^{\nu}\\&\qquad\qquad+\sum_{\nu=0}^{m+1}|\beta_\nu t_1t_2+\beta_{\nu-1}(t_1-t_2)-\beta_{\nu-2}|t_{1}^{\nu}\\&\qquad\qquad+\sum_{\nu=m+2}^{n}|\beta_\nu t_1t_2+\beta_{\nu-1}(t_1-t_2)-\beta_{\nu-2}|t_{1}^{\nu}\\&= 2(\alpha_{k+1}t_2+\alpha_k)t_1^{k+2}-(\alpha_nt_2+\alpha_{n-1})t_1^{n+1}\\&\qquad+2(\beta_{m+1}t_2+\beta_m)t_1^{m+2}-(\beta_nt_2+\beta_{n-1})t_1^{n+1}.
\end{align*}
Using this in \eqref{p31}, we obtain
\begin{align*}
|G(z)|\geq& |z|^{n+1}\bigg\{\left|a_nz+\alpha_{n-1}-(t_1-t_2)\alpha_n\right|-(\beta_{n-1}-(t_1-t_2)\beta_n)\\&-2(\alpha_{k+1}t_2+\alpha_k)\dfrac{1}{t_1^{n-k-1}}+(\alpha_nt_2+\alpha_{n-1})\\&-2(\beta_{m+1}t_2+\beta_m)\dfrac{1}{t_1^{n-m-1}}+(\beta_nt_2+\beta_{n-1}) \bigg\}\\=&|z|^{n+1}\bigg\{\left|a_nz+\alpha_{n-1}-(t_1-t_2)\alpha_n\right|+t_1\beta_n\\&-2(\alpha_{k+1}t_2+\alpha_k)\dfrac{1}{t_1^{n-k-1}}+(\alpha_nt_2+\alpha_{n-1})\\&\qquad\qquad-2(\beta_{m+1}t_2+\beta_m)\dfrac{1}{t_1^{n-m-1}} \bigg\}>0,
\end{align*}
if
\begin{align*}\nonumber
\big|a_nz+&\alpha_{n-1}-(t_1-t_2)\alpha_n\big|\\\nonumber&>2(\alpha_{k+1}t_2+\alpha_k)\dfrac{t_1^{k+1}}{t_1^{n}}+2(\beta_{m+1}t_2+\beta_m)\dfrac{t_1^{m+1}}{t_1^{n}}\\&\qquad\qquad\qquad-(t_2\alpha_n+t_1\beta_n+\alpha_{n-1}).
\end{align*}
Hence all the zeros of $G(z)$ whose modulus is greater than $t_1$ lie in the circle 
\begin{align}\nonumber\label{p11}
\bigg|z+&\dfrac{\alpha_{n-1}-(t_1-t_2)\alpha_n}{a_n}\bigg|\\\nonumber&\leq 2(\alpha_{k+1}t_2+\alpha_k)\dfrac{t_1^{k+1}}{|a_n|t_1^{n}}+2(\beta_{m+1}t_2+\beta_m)\dfrac{t_1^{m+1}}{|a_n|t_1^{n}}\\&\qquad\qquad\qquad-\dfrac{t_2\alpha_n+t_1\beta_n+\alpha_{n-1}}{|a_n|}.
\end{align}
Now, we show that all the zeros of $G(z)$ whose modulus is less than equal to $t_1$ also lie in the circle defined by \eqref{t2e}. Let $|z|\leq t_1,$ then we have 
\begin{align*}\nonumber
\big|a_nz+\alpha_{n-1}-(t_1-t_2)\alpha_n\big|&\leq |a_n|t_1+|\alpha_{n-1}-(t_1-t_2)\alpha_n|\\\nonumber&\leq t_1\alpha_n+t_1\beta_n+\alpha_{n-1}-(t_1-t_2)\alpha_n\\\nonumber&=2t_1\beta_n+2(t_2\alpha_n+\alpha_{n-1})\\&\qquad-(t_2\alpha_n+2t_1\beta_n+\alpha_{n-1}).
\end{align*}
By hypothesis,
\begin{align*}\nonumber
\sum_{\nu=k+2}^{n}\dfrac{\alpha_\nu t_1t_2+\alpha_{\nu-1}(t_1-t_2)-\alpha_{\nu-2}}{t_1^{n-\nu+1}}\leq 0,\qquad 0\leq k\leq n-1,
\end{align*}
this gives
\begin{align}\label{p7}
 2\left(t_2\alpha_n+\alpha_{n-1}\right)\leq \dfrac{2t_2\alpha_{k+1}+2\alpha_k}{t_1^{n-k-1}}
\end{align}
 Similarly for $0\leq m\leq n-1,$
\begin{align}\label{p8}
2\left(t_2\beta_n+\beta_{n-1}\right)\leq \dfrac{2t_2\beta_{m+1}+2\beta_m}{t_1^{n-m-1}} .
\end{align}
Also, we have
\begin{align}\label{p9}
 t_1\beta_n\leq t_2\beta_n+\beta_{n-1}
\end{align}
Combining  \eqref{p8} and \eqref{p9}, we obtain
\begin{align}\label{p10}
2 t_1\beta_n\leq \dfrac{2t_2\beta_{m+1}+2\beta_m}{t_1^{n-m-1}}
\end{align}
Using \eqref{p7} and \eqref{p10} in \eqref{p11}, we have
\begin{align*}
\big|a_nz+&\alpha_{n-1}-(t_1-t_2)\alpha_n\big|\\\nonumber&\leq 2(\alpha_{k+1}t_2+\alpha_k)\dfrac{t_1^{k+1}}{t_1^{n}}+2(\beta_{m+1}t_2+\beta_m)\dfrac{t_1^{m+1}}{t_1^{n}}\\&\qquad\qquad\qquad-(t_2\alpha_n+t_1\beta_n+\alpha_{n-1}).
\end{align*}
This shows that all the zeros of $G(z)$ whose modulus is less than or equal to $t_1$ also satisfy the inequality \eqref{t2e}. Thus we conclude that all the zeros of $G(z)$ hence that of $P(z)$ lie in the circle defined by \eqref{t2e}.  This completes the proof of Theorem 1.10.
\end{proof}


\begin{thebibliography}{99}
{\smaller
\bibitem{asv79} N. Anderson, E. B. Saff, R. S. Verga, On the Enestr\"{o}m-Kakeya theorem and its sharpness, \textit{Linear Algebra and Applications} \textbf{28} (1979), 5-16.
\bibitem{asv81} N. Anderson, E. B. Saff, R. S. Verga, An extension of the Enestr\"{o}m-Kakeya theorem and its sharpness, \textit{SIAM J. Math. Anal.,} \textbf{12} (1981), 10-22.
\bibitem{am} A. Aziz and Q. G. Mohammad, Zero free regions for polynomials and some generalizations of Enestr\"{o}m-Kakeya theorem, \textit{Canad. Math. Bull.,}, \textbf{2} (1984) 265-272.
\bibitem{cs} G. T. Cargo and O. Shisha,  Zeros of polynomials and fractional order differences of their coefficients, \textit{J. Math. Anal. Appl.,} \textbf{7} (1963), 176-182.
\bibitem{e} G. Enestr\"{o}m, Ramarque sur un th\'{e}or\`{e}me relatif aux recines de l'equation $a_nx^n+\cdots+a_0=0$ o\`{u} tous les coefficients sont r\'{e}els et positifs, \textit{T\^{o}hoku Math. J.,} \textbf{18} (1920), 34-36; translation of a Swedish article in \textit{Ofversigt of Konogl. Vertenskaps  Akademiens F\"{o}rhandlingar,} \textbf{50} (1893), 405-415.
 \bibitem{gr} N. K. Govil and Q. I. Rahman, On the Enestr\"{o}m-Kakeya theorem , \textit{Tohuku Math. J.}, \textbf{20} (1968) 126-136
 \bibitem{sk} S. Kakeya, On the limits of the roots of an algebraic equation with positive coefficients, \textit{Tohoku Math. J.,} \textbf{2} (1912-1913), 140-142.
\bibitem{km} M. Kova\v{c}evi\'{c} and I. Milovanovi\'{c} , On a generalization of the Enestr\"{o}m-Kakeya theorem, \textit{Pure Math. and Applic. Ser. A}  \textbf{3} (1992), 543-47.
\bibitem{k} P. V. Krishnaih, On the Enestr\"{o}m-Kakeya theorem, \textit{J. London Math. Soc.,} \textbf{30} (1955) 314-319.
\bibitem{marden} M. Marden, \textit{Geometry of Polynomials}, Math. Surveys No. 3, Amer. Math. Soc. Providence R. I. 1966.
\bibitem{mmr}G. V. Milovanovic, D. S. Mitrinovic and Th. M. Rassias, \textit{Topics in Polynomials: Extremal Properties, Inequalities}, Zeros, World scientific Publishing Co., Singapore, (1994).
\bibitem{rsm} N. A. Rather, Shakeel A. Simnani, M. I. Mir, On the Enestr\"{o}m-Kakeya theorem, \emph{Int. J. Pure and Appl. Math.,} \textbf{41} (2007) 807-815.
\bibitem{rsuh} N. A. Rather and Suhail Gulzar, On the Enestr\"{o}m-Kakeya theorem, \emph{Acta Math. Univ. Comenianae,} \textbf{83} (2014) 291-302.
\label{pgCAbdio}
}
\end{thebibliography}
\end{document}